\documentclass[12pt]{amsart}
\usepackage{amsmath, amscd, amssymb, latexsym}

\theoremstyle{plain}
\newtheorem{thm}{Theorem}[section]
\newtheorem{lem}[thm]{Lemma}
\newtheorem{cor}[thm]{Corollary}
\newtheorem{pro}[thm]{Proposition}

\theoremstyle{definition}
\newtheorem{df}[thm]{Definition}

\newtheorem{rem}[thm]{Remark}
\newtheorem{ex}[thm]{Example}

\def\Om{\Omega}
\def\La{\Lambda}

\def\ov{\overline}

\def\al{\alpha}
\def\ga{\gamma}

\def\na{\nabla}

\def\wt{\widetilde}
\def\la{\lambda}
\def\ni{\noindent}
\def\nt{\noindent}

\def\we{\wedge}

\def\un{\underline}

\def\pa{\partial}

\def\mc{\mathcal}
\def\mf{\mathfrak}
\def\i{{\bf i}}
\def\mr{\mathrm}
\def\R{{\bf R}}
\def\C{{\bf C}}
\def\Z{{\bf Z}}

\begin{document}

\title{OBJECTIVE B-FIELDS AND A HITCHIN-KOBAYASHI CORRESPONDENCE}

\author{SHUGUANG WANG}
\address{Department of Mathematics
\\University of Missouri\\Columbia,
MO 65211,USA}

\email{wangs@missouri.edu}

\begin{abstract}
 A simple trick invoking objective B-fields is employed to refine
the concept of characteristic classes for twisted bundles.
Then the objective stability and objective Einstein metrics are introduced and a new
Hitchin-Kobayashi correspondence is established between them. As an application
the SO(3)-instanton moduli space is proved to be always orientable.
\end{abstract}

\maketitle

\nt{\em Keywords}: Objective Chern classes/Hermitian-Einstein metric/stable bundle;
Hitchin-Kobayashi correspondence; orientable SO(3)-instanton moduli space.\\


\vspace{5mm}

\centerline{\sc 0. Introduction}
\vspace{3mm}

Recently there has been intensive literature on gerbes starting from
 Brylinski's book \cite{b}, in which the Chern-Weil theory of
gerbes was carried out. The geometry of abelian gerbes was further refined
and clarified by Murray \cite{m} and Chatterjee \cite{c}, Hitchin \cite{h}
 from two different view-points. These objects have found
applications in Physics for the descriptions of the so-called B-fields and
  the twisted K-theory according to \cite{a, bm, lu,pry} among many others.
  A real version of gerbes is presented in \cite{w2}.

In one aspect gerbes are utilized to define twisted bundles where the usual
cocycle condition
 fails to hold.  The concept of twisted principal bundles was  explicitly
 proposed in \cite{ma,mp} where their holonomy is studied, while
 the authors in \cite{bs,ms} introduced bundle gerbe modules in a different
setting for the generalization of the index theory. In the context of
algebraic geometry
these are Azumaya bundles  via the description of
gerbes in terms of Brauer groups, see \cite{mms1,y,ca} for example.

The first main purpose in our paper is to propose a new concept of
Chern classes for twisted vector bundles as a refinement of
existing definitions from \cite{bs, ca}. In fact the earlier
definitions  are  not completely satisfactory in the sense that
the one give in \cite{bs} is not quite topological because of the
dependence on different curvings, while the definition in
\cite{ca} involves an un-specified choice of a second twisted
bundle. Our approach can be viewed as a balance of the two. To
explain the key idea behind, recognizing that  the underlying
gerbe $\mc{L}$ of a rank $m$-twisted vector bundle $E$ is an
$m$-torsion, we select a trivialization $L$ of $\mc{L}^m$. To remove the
dependence on the B-fields $B$ (namely curvings), we restrict to a special
$B$ so that $mB$ is objective with respect to $L$.
Given a twisted connection $D$ on $E$, twisted over a gerbe connection $\mc{A}$,
we use their curvatures and $B$ together to define closed global Chern forms
$\phi_k(E,L;\mc{A},B,D)$.\\

\ni{\bf Theorem 0.1.} {\em The  $L$-objective Chern classes of $E$,  defined as
 $c_k^L(E)=[\phi_k(E,L;\mc{A},B,D)]\in H^{2k}(X,\R)$,
are topological invariants of $E$. In other words, $c_k^L(E)$
are independent of the choices made for $D,\mc{A},B$.}\\

\nt(The usage of the word ``objective'' throughout the paper is inspired by \cite{c}.)
It seems likely that Theorem 0.1 will provide an alternative starting point to
the current research \cite{bm,bs,ms,y} concerning a host of index theorems and
twisted K-theory.

With objective Chern classes and properly defined twisted subbundles, we
further introduce
objective stability. Moreover by fixing a metric on $L$, we are able to
characterize objective Einstein metrics so that we can formulate and generalize
the Hitchin-Kobayashi correspondence as follows. (A possible statement as such
 was  alluded in \cite{dr}.) The original correspondence for untwisted bundles
was a corner-stone result that  had been established in \cite{l,d1,uy,d2}.\\

\ni{\bf Theorem 0.2.} {\em Suppose $E$ is a  rank-$m$ indecomposable
holomorphic bundle twisted over a gerbe $\mc{L}$ and $L$ is a holomorphic
trivialization of $\mc{L}^m$. Then $E$ is $L$-objective stable iff $E$
admits an $H$-objective
Einstein metric for some Hermitian metric $H$ on $L$.}\\

The theorem is new only if the base manifold has complex dimension
at least $2$. For the case of a complex curve, any gerbe $\mc{L}$
must be trivial hence no twisting occurs. But one might still ask
about stable bundles and Yang-Mills connections on non-orientable
Riemann surfaces, compare \cite{w1} for example.

Our consideration of trivializations $L$ of $\mc{L}^m$ is rather natural and
motivated in part by
the lifting of $\mr{SO}(3)$-bundles. More precisely when an  $\mr{SO}(3)$-bundle
$S$ is $\mr{spin}^c$, lifting its structure group from
$\mr{SO}(3)$ to $\mr{U}(2)$ one has the spinor bundle $E$ with determinant line
bundle $\det E$. In general if $S$ is not $\mr{spin}^c$, one can only have a twisted
$\mr{U}(2)$-bundle $E$ over a non-trivial gerbe $\mc{L}$. Then a trivialization $L$
of $\mc{L}^2$ is a replacement of the determinant.
 In this regard,  we get the following result  as an application of our theory:\\

\ni{\bf Proposition 0.3.} {\em The  moduli space of  anti-self dual
connections on any $\mr{SO}(3)$-bundle $S$ over a 4-manifold is orientable,
whenever it is smooth.}\\

\nt This is a small generalization of Donaldson's theorem \cite{d3} where $S$ is assumed
to be $\mr{spin}^c$. Donaldson's theorem has proved to be quite crucial
for various gauge theoretic applications.

Here is a brief description of the paper. After reviewing and setting up notations
for gerbes and so on in Section \ref{ger}, we develop the objective Chern-Weil theory
and prove Theorem 0.1 in Section \ref{occt}. The fundamental concepts of
objective subbundles, Einstein metrics
 and stability are laid out in Sections \ref{oss} and \ref{thea}. Then we devote
 Sections \ref{ghk}, \ref{oth} to the proofs of the two directions in Theorem 0.2.
Finally Section \ref{aap} contains the proof of Proposition 0.3.\\

\nt{\bf Acknowledgments.} We would like to thank Yongbin Ruan for several contributions
and discussions that proved to be quite helpful. The early portion of the
paper was written while the author was visiting IHES in the summer of 2008.

\section{Review of gerbes, B-fields and twisted bundles}\label{ger}

All vector bundles, differential forms will be defined on the
complex field $\C$ unless otherwise indicated.

To set up notations, first recall the definition of gerbes from
\cite{c,h}. Let $X$ be a smooth manifold $X$ bearing an  open
cover $\{U_i\}$. A {\em gerbe} $\mc{L}=\{L_{ij}\}$  consists of
line bundles $L_{ij}\to U_i\cap U_j$ with  given isomorphisms
\begin{equation}\label{m2}
 L_{ij}=L_{ji}^{-1},\;\;\;L_{ij}\otimes L_{jk}= L_{ik}
 \end{equation}
  on their common domains of definition.
Using a good or  refinement cover if necessary, we can assume
  $L_{ij}$ are trivial bundles. Choose
trivializations $\xi_{ij}$ for $L_{ij}$. They should be compatible
with (\ref{m2}) in that
\begin{equation}\label{m3}
\xi_{ij}=\xi_{ji}^{-1},\;\;\;\xi_{ij}\otimes\xi_{jk}=z_{ijk}\xi_{ik}
\end{equation}
where $\un{z}=\{z_{ijk}\}$ forms a $\check{\mbox{C}}$ech 2-cycle
of the sheaf $\un{\C}^*$ of nowhere vanishing complex functions on
$X$. One imposes that $\un{z}$ be co-closed: $\delta \un{z}=1$.

Thus the {\em gerbe class} $[\mc{L}]:=[\un{z}]\in\check{H}^2(X,
\un{\C}^*)$ is defined. Two gerbes are {\em isomorphic} if they
become the same on a common refinement. This is equivalent to
their gerbe classes being the same. As in the case of bundles, it
is often convenient  not to distinguish a gerbe $\mc{L}$ from its
class $[\mc{L}]$. One can
 view $[\mc{L}]\in {H}^3(X, {\bf Z})$ according to the natural
isomorphism $\check{H}^2(X, \un{\C}^*)\to {H}^3(X, {\bf Z})$,
which comes from the short exact sequence ${\bf
Z}\to\un{\C}\to \un{\C}^*$. Gerbes form an Abelian group under
the tensor product.\\

From (\ref{m2}), one has a trivialization $s_{ijk}$ of
$L_{ijk}:=L_{ij}\otimes L_{jk}\otimes L_{ki}$ as part of the gerbe
definition. In \cite{c,h}, a gerbe connection $\{A_{ij}, B_i\}$ on
$\mc{L}$ consists of two parts: connections $A_{ij}$ on $L_{ij}$
and complex 2-forms $B_i$ on $U_i$, such that 1) $s_{ijk}$ is a
covariant constant trivialization of $L_{ijk}$ under the induced
product connection using $A_{ij}$, and 2) $B_j-B_i=F_{A_{ij}}$ on
$U_i\cap U_j$ for the curvature of $A_{ij}$. However in this paper
we will separate them and call $\mc{A}=\{A_{ij}\}$ alone  a {\em
gerbe connection} and give a more prominent role to $B=\{B_i\}$,
which we call a ($\mc{A}$-compatible) {\em B-field} from its
physics interpretation.  The closed local 3-form $\{d B_i\}$ can
be patched together to yield a global form $G$, called the {\em
gerbe curvature} of the pair $\mc{A}, B$. The class $[G]\in
H^3(X,\R)$ is the real image of $[\mc{L}]\in H^3(X,\Z)$, hence
independent of the choice of $\mc{A}, B$. Moreover every
representative of the real image can be realized by some B-field
and connection. Here is a collection of useful facts to recall:

\begin{pro}\label{tor}
 The following are equivalent.

{\em 1)}  $\mc{L}$ is a torsion, i.e. $m[\mc{L}]=0\in H^3(X,\Z)$
for some integer $m$.

{\em 2)} Under a suitable cover of $X$, $[\mc{L}]$ can be
represented by a $\check{\mbox{C}}$ech 2-cycle
$\un{z}=\{z_{ijk}\}$  such that $z_{ijk}^m=1$ for some fixed $m$.
(Hence $\un{z}$ must be locally constant.)

{\em 3)} The real gerbe class of $\mc{L}$ is trivial, namely the
image $[\mc{L}]=0$ in $H^3(X,\R)$.

{\em 4)} $\mc{L}$ is flat, i.e. there exists a  gerbe connection
$\mc{A}$ with a compatible B-field $B$ such that
 the curvature 3-form $G=0$.

\end{pro}

Next suppose $\mc{L}$ is trivial, i.e. its class $[{\mc L}]=0\in
\check{H}^2(X, \un{\C}^*)$. Then $\un{z}$ is a coboundary of a
$\check{\mbox{C}}$ech 1-cocycle, which leads to a {\em
trivialization} $L=\{L_i\}$ of $\mc{L}$, namely line bundles $L_i$
on $U_i$ such that $L_j=L_i\otimes L_{ij}$ when restricted to
$U_i\cap U_j$. Given a second trivialization $L'=\{L'_i\}$, we
have a global line bundle $\ell$ on $X$, called the {\em
difference bundle} of $L$ and $L'$. We will denote $\ell=L\ominus
L'$. A {\em trivialization} of a gerbe connection $\mc{A}$ is a
family of connections $A=\{A_i\}$ on $\{L_i\}$ subject to an
analogous condition. Since the 1-form sheaf is a fine sheaf,
trivializations always exist for any gerbe connection but not
unique. Given a second trivializing connection $A'$ on $L'$, we
have a {\em difference connection} $A\ominus A'$ on the bundle
$\ell$.

Given a trivializing connection $A$ of $\mc{A}$, the collection of
local curvature forms $\{F_{A_i}\}$ is a B-field compatible with
$\mc{A}$; but not every compatible B-field  has this form (i.e.
from some trivializing connection  of $\mc{A}$). Given a
compatible B-field $B$ and a trivializing connection $A$, the
difference $\epsilon:=B_i-F_{A_i}$ is a global 2-form on $X$. This
{\em error form} $\epsilon$ is essential in defining the gerbe
holonomy of $\mc{A}$. Following \cite{c} we call $B$  {\em objective}
if $\epsilon=0$ for some $A$. A necessary condition is
that $B_i$ be closed. Suppose so, then the error form $\epsilon$
with any trivializing connection  is also closed hence gives rise
to a class $[\epsilon]\in H^2(X,\R)$. In this language,  objective
B-fields are exactly those closed B-fields such that the error
class $[\epsilon]$ lives in the lattice $H^2(X,\Z)\subset
H^2(X,\R)$. Here of course a gerbe connection $\mc{A}$ has been
fixed throughout.

Now is a good time to set out the following guide for our
notations in the paper:

\begin{rem} We will continue using script letters such as $\mc{L},\mc{A}$
for things defined on intersections of open sets $\{U_i\}$ and
reserve standard block letters $L, B, A$ etc for objects defined
on open sets themselves.
\end{rem}

A trivialization $L=\{L_i\}$ of $\mc{L}$ can also be viewed as a
twisted line bundle over $\mc{L}$. In general, given an arbitrary
gerbe $\mc{L}$, an {\em $\mc{L}$-twisted vector bundle}
$E=\{E_i\}$ on $X$ consists of a collection of local bundles of
the same rank such that $E_j=E_i\otimes L_{ij}$ on $U_i\cap U_j$.
We will also say that $E$ is over the {\em twisting gerbe} ${\mc
L}$. (Twisted principal bundles are first defined in \cite{ma}.)
 In a similar spirit, one defines a {\em twisted
connection} $D=\{D_i\}$ on $E$ over a gerbe connection $\mc{A}$.
We will adopt the following convenient notations to indicate the
twisting gerbe $\mc{L}$ and the twisting connection $\mc{A}$:
$$E\prec{\mc L},D\prec\mc{A}.$$

Twisted bundles/connections are not necessarily mysterious: simply
put, these are so defined that their projections are just regular
fiber bundles and connections respectively.  For any 1-dimensional
vector space $V$, $V\otimes V^*$ carries a canonical basis. It
follows that the local endomorphism bundles $\{\mr{End}(E_i)\}$
naturally fit together to produce a global vector bundle on $X$,
in view of $\mr{End}(E_i)=E_i\otimes E^*_i$. We denote it by
$\mr{End}(E)$. Note that the wedge product twisted bundle is over
$\mc{L}^r$:
$$\we^r E=\{\we^r E_i\}\prec\mc{L}^r=\{\otimes^rL_{ij}\}.$$
In particular $\det E$ gives rise to a trivialization of $\mc{L}^m$, where
$m=\mr{rank}E$. In other words, the underlying gerbe ${\mc L}$ of
$E$ must be an $m$-torsion.

Other operations can be introduced as well: one has the dual
$E^*\prec\mc{L}^{-1}$. If $E'\prec\mc{L}'$ then
$E\otimes E'\prec\mc{L}\otimes\mc{L}'$. But for $E\oplus E'$,
one must impose $\mc{L}=\mc{L}'$; then $E\oplus E'\prec\mc{L}$.

\section{Objective Chern classes for twisted bundles}\label{occt}

We now develop a new version of Chern-Weil theory for rank
$m$-twisted vector bundles. Fix  a gerbe $\mc{L}$ with connection
$\mc{A}$ and B-field $B$. Assume $\mc{L}$ is an $m$-torsion. By
Proposition \ref{tor}, we can choose $B$ to be flat so that all
$B_i$ are closed 2-forms. For a twisted rank $m$-vector bundle
$E\prec\mc{L}$ with connection $D\prec\mc{A}$, it is easy to check
that {\em B-twisted curvature} of $D$:
$$\wt{F}_D = F_{D_i}-B_i I$$
is a global section of $\Om^2(\mr{End}(E))$. In order to get the
expected topological results we will need to place a crucial
restriction on $B$.

\begin{thm}\label{chern}Consider an $m$-torsion
 gerbe $\mc{L}$ with connection $\mc{A}$ and B-filed $B$. Let
$L=\{L_i\}$ be a trivialization of $\mc{L}^m$ and $A=\{A_i\}$ be
that of $\mc{A}^m$. Assume the B-field $mB$ on $\mc{L}^m$ to be
objective: $mB=\{F_{A_i}\}$ (one may call $B$ relative objective).

Let $E$ be a rank-$m$ twisted vector bundle over $\mc{L}$ and
$D=\{D_i\}$  a twisted connection over $\mc{A}$. Introduce the
total Chern form
$$\phi(E,L;\mc{A},B,D)=\det(I+\frac{i}{2\pi}\wt{F}_D).$$

{\em 1)} Each $\phi_k(E,L;\mc{A},B,D)$ is a closed $2k$ form on $X$.

{\em 2)} The  class $[\phi_k]\in H^{2k}(X,\R)$ is independent of
the choices of $\mc{A},B, D$.

{\em 3)} Suppose $L'=\{L_i'\}$ is a second trivialization of
$\mc{L}^m$ and $A'=\{A'_i\}$ is a trivialization of $\mc{A}^m$ on
$L'$. Let $B'$ be the B-field on $\mc{L}$ such that
$mB'=\{F_{A'_i}\}$.
 Then the following holds:
 $$[\phi_k(E,L';\mc{A},B',D)]=\sum^k_{i=0}\frac{C(m-i,k-i)}{m^{k-i}}
 [\phi_i(E,L;\mc{A},B,D)]c_1(\ell)^{k-i}$$
where $\ell=L\ominus L'$ is the difference line bundle of $L, L'$.
\end{thm}

\nt{\it Proof.} 1) Following the standard Chern-Weil theory (for
example from Chapter III of \cite{wl}), it suffices to show a
Bianchi type identity still holds:
\begin{equation}\label{bi}
d\wt{F}_D=[\wt{F}_D,\wt{\theta}] \end{equation}
 for some
$\wt{\theta}\in\Omega^1(\mr{End}(E))$ on any small open set
$V\subset U_i$. To prove the identity, if $V$ is small enough,
there is a line bundle $K\to V$ such that $L_i=K^m$. Set $A_i=G^m$
for some unique connection $G$ on $K$. Then
$B_i=\frac{1}{m}F_{A_i}=F_G$. It follows that
$$\wt{F}_D=F_{D_i}-F_GI=F_{\wt{D}}$$
for the curvature of the tensor product connection
$\wt{D}=D\otimes G^{-1}$ on $E_i\otimes K^{-1}$. Now the standard
Bianchi identity for $\wt{D}$ yields (\ref{bi}), where
$\wt{\theta}$ is the connection matrix of $\wt{D}$ under a local
frame.

2) Note that $\mc{A}$ and $B$ are both determined uniquely by $A$
(although $\mc{L}$ is not so by $L$). It is enough to show more
generally that $[\phi_k]$ is independent of the choices of $D, A$.
For a second pair of choices $D'=\{D'_i\}, A'=\{A'_i\}$, consider
the one-parameter families $D(t)=\{tD_i+(1-t)D'_i\},
A(t)=\{tA_i+(1-t)A'_i\}$. Then $A(t)$ determines  families of
gerbe connections $\mc{A}(t)$ and objective B-fields $B(t)$ on
$\mc{L}$. Moreover $D(t)$ is a twisted connection over $\mc{A}(t)$
for each $t$.

Set the family $\wt{F}(t)=F_{D(t)}-B(t)I$. On a small enough open
set $V\subset U_i$, let $L_i=K^m$ as above and $G(t)$ be the
family of connections on $K$ determined by $A_i(t)$. Then
$\wt{F}(t)=F_{\wt{D}(t)}$, the curvature of the tensor product
connections $\wt{D}(t)=D(t)\otimes G(t)^{-1}$ on $E_i\otimes
K^{-1}$. Thus
$$\wt{F}(t)=d\wt{\theta}(t)+\wt{\theta}(t)\wedge\wt{\theta}(t)$$
for the connection matrix $\wt{\theta}(t)$ of $\wt{D}(t)$ under a
frame. Consequently
$$\dot{\wt{F}}=d\dot{\wt{\theta}}+[\dot{\wt{\theta}},\wt{\theta}].$$
As in the standard Chern-Weil theory \cite{wl}, the above formula
together with the Bianchi formula (\ref{bi}) shows that
$$\phi_k(E,L;\mc{A}',B',D')-\phi_k(E,L;\mc{A},B,D)$$
is an exact form. Hence
$[\phi_k(E,L;\mc{A}',B',D')]=[\phi_k(E,L;\mc{A},B,D)]$.

3)  Here it is important to be able to keep the same gerbe
connection $\mc{A}$ on $\mc{L}$ and hence the same twisted
connection $D$ on $E$. In other words a different trivialization
$L'$ will only impact on the B-field  via $A'$. Now
$$F_D-B'I=F_D-BI+(B-B')I=F_D-BI+\frac{1}{m}F_\alpha I,$$ where $\alpha$ is
the difference connection of $A,A'$ on $\ell$. Under a frame,
write $F_D-BI$ as a matrix $(F^i_j)$ of 2-forms. Then $F_D-B'I$ is
given by the matrix $(F^i_j+\frac{1}{m}\delta^i_j F_\al)$.
Applying the classical formula to $\det(I+F_D-BI)$, one has
$$\phi_k(E,L;\mc{A},B,D)=\frac{1}{(2\pi i)^k k!}\sum
\delta^{j_1\cdots j_k}_{i_1\cdots i_k}
F^{i_1}_{j_1}\wedge\cdots\wedge F^{i_k}_{j_k}.$$
 Doing the same for $\det(I+F_D-B'I)$ will give
$$\phi_k(E,L';\mc{A},B',D)=\frac{1}{(2\pi i)^k k!}\sum
\delta^{j_1\cdots j_k}_{i_1\cdots i_k}
(F^{i_1}_{j_1}+\frac{1}{m}\delta^{i_1}_{j_1}
F_\al)\wedge\cdots\wedge
(F^{i_k}_{j_k}+\frac{1}{m}\delta^{i_k}_{j_k} F_\al).$$
 These formulas together will give us the desired result for their
 classes after noting $c_1(\ell)=[F_\al]$.
  q.e.d.

In particular if $L'$ is flat equivalent to $L$, namely if $\ell$
is flat, then $c_1(\ell)=0\in H^2(X,\R)$ and the class $[\phi_k]$
remains the same. This is the case if $L'$ is isomorphic to $L$
(i.e. $\ell$ is trivial). Thus it makes sense to introduce:

\begin{df}\label{tc}
Suppose $E\prec\mc{L}$ is a rank-$m$ twisted bundle and
$L\prec\mc{L}^m$ is a trivialization. The
 {\em $L$-objective Chern class} of $E$ is defined as
follows: $$c_k^L(E)=[\phi_k(E,L;\mc{A},B,D)]\in H^{2k}(X,\R),$$
where $\mc{A}$ is any gerbe connection on $\mc{L}$, $D$ is any
twisted connection on $E$ and $B$ is a B-field compatible with
$\mc{A}$ such that $mB$ is objective.
\end{df}

\begin{rem}
Alternatively one could define a twisted bundle as a {\em pair}
$(E,L)$ where $E\prec\mc{L}$ and $L\prec\mc{L}^m$ in our
notations. Then one could define the Chern classes for $(E,L)$
without the cumbersome prefix ``$L$-objective''. For us the
advantage of separating $L$ from $E$ is that we can illustrate
better the dependence on $L$ as shown in the next corollary.
\end{rem}

\begin{cor}\label{c1}
 Suppose $E\prec\mc{L}$ and $L,L'\prec\mc{L}^m$. Then
$$c_k^{L'}(E)=\sum^k_{i=0}\frac{C(m-i,k-i)}{m^{k-i}}
c_i^{L}(E)c_1(\ell)^{k-i}$$ in terms of the difference line bundle
$\ell=L\ominus L'$. In particular
$$c^{L'}_1(E)=c^L_1(E)+ c_1(\ell).$$
Recall $\wedge^r E\prec\mc{L}^r$ for any integer $r\leq m$.
Setting $w=C(m-1,r-1)$, one has also
$$c_1^{L^w}(\wedge^r E)=wc_1^L(E).$$
In particular $c_1^{L}(\mr{det} E)=c_1^L(E)$ as one might expect.
\end{cor}

\nt{\it Proof.} The first formula is a translation of Part 3) of
Theorem \ref{chern}. To show the third formula, note that
$\mr{rank}\wedge^r E=C(m,r)$ and
$$rC(m,r)=mC(m-1,r-1)=mw.$$
Hence $L^w$ is a trivialization of $\mc{L}^{rC(m,r)}$, which
carries the B-field $C(m,r)rB$. Let $\wedge D$ denote the induced
twisted connection on $\wedge^r E$. One has
$$\begin{array}{ll}
c_1^{L^w}(\wedge^r E)&=\mr{tr} F_{\wedge D}-C(m,r)rB\\
&=w\mr{tr} F_D-mwB=w(\mr{tr} F_D-mB)\\
& =wc^L_1(E)
\end{array}$$
delivering us the desired formula.  q.e.d.

\begin{rem}\label{ch}
One could in particular use the special trivialization  $\mr{det}E=
\wedge^m E\prec\mc{L}^m$
 to introduce $c_k(E):=c^{\mr{det}E}_k(E)$. For any
other trivialization $L\prec\mc{L}^m$, the previous corollary then
yields
$$c_k^{L}(E)=\sum^k_{i=0}\frac{C(m-i,k-i)}{m^{k-i}}
c_i(E)c_1(\ell)^{k-i}$$ using the difference bundle
$\ell=\mr{det}E\ominus L$. Note $c_1(E)=0$ and the formula above
implies that $c_1^L(E)=c_1(\ell)$ {\em is always an integer
class}.
\end{rem}

The following example from gauge theory has been a useful guide
for us and will also show why it is important to expand $c_k(E)$
to $c^L_k(E)$ with the incorporation of  trivializations other
than $\mr{det}E$.

 \begin{ex}\label{e1}
  Let ${Q}$ be an $\mr{SO}(3)$-vector bundle on $X$.
 Suppose $Q$ is $\mr{spin}^c$ so that $w_2(Q)$ has
 integral lifts in $H^2(X,\Z)$. Fix such a lift $L$ viewed
 as a line bundle as well as a connection $A$ on it.
 Then $Q$ together with  connection $\na$  will lift to a unique
 $\mr{U}(2)$-vector bundle $\wt{Q}$ with connection $\wt{\na}$
 such that  $\mr{det}(\wt{Q})=L,\mr{det}{\wt{\na}}=A$.  We can
 interpret $\wt{Q},\wt{\na}$ and especially the Chern classes
 $c_k(\wt{Q})$
in terms of our twisted vector bundle theory. For this purpose
 choose any open cover $\{U_i\}$ of $X$ such that $L|_{U_i}$ has
 a square root $K_i:=\sqrt{L}$ on $U_i$. Thus $K=\{K_i\}$ is
 a trivialization
of the gerbe $\mc{L}=\{L_{ij}\}$ where $L_{ij}=K_j\otimes K^*_i$
by definition. (The gerbe $\mc{L}$, first observed in \cite{we},
can be formally viewed as being given by $\sqrt{L}$, namely
$\mc{L}$ is the obstruction to the existence of a global squared
root of $L$. Similarly there is a gerbe given by the $n$-th root
$\sqrt[n]{L}$. More on real gerbes can be found in \cite{w1}.)
Moreover $A$ leads to a unique twisted connection
$\wt{A}=\{\wt{A}_i\}\prec\mc{A}$ on $K$, where $\mc{A}=\{A_{ij}\}$
and $A_{ij}$ is the tensor product connection $\wt{A}_j\otimes
\wt{A}_i^*$ on $L_{ij}$. Now on each $U_i$, lift $\wt{Q}$ to an
$\mr{SU}(2)$ vector bundle $E_i$ so that $E=\{E_i\}$ is a twisted
bundle over $\mc{L}$. Then $\na$ lifts to a unique twisted
connection $D=\{D_i\}\prec\mc{A}$ on $E$. Now one can check easily
that
$$\wt{Q}=E_i\otimes K_i, \wt{\na}=D_i\otimes\wt{A}_i,
c_k(\wt{Q})=c^L_k(E),$$
 where $L=K^2\prec\mc{L}^2$ is viewed as a
trivialization of $\mc{L}^2$. (Here $K_i^2$ is global since
$\mc{L}^2$ is trivial as a $\check{\mr{C}}$ech cycle not just as a
$\check{\mr{C}}$ech class.)
\end{ex}

When $Q$ is not $\mr{spin}^c$, the line bundle $L$ will not exist any
more. Nonetheless, $Q$ still  lifts locally to a twisted $\mr{SU}(2)$
bundle $E=\{E_i\}$ over some gerbe $\mc{L}$. (The non-trivial
gerbe class $[\mc{L}]$ equals $W_3(Q)\in H^3(X,\Z)$.) The
{\bf essence} of our theory above is {\em to replace $L$ with a trivializing
twisted bundle of $\mc{L}^2$, and such a trivialization always
exits!}  More  generally, the lesson for our entire paper is this:
Given twisted data such as $E, D$
that are defined on local open sets only, one should couple them with  {\bf
objective} gerbe data in order to get global data on $X$.

\begin{rem}\label{ch3} We now compare our Definition \ref{tc} with several
existing definitions of Chern classes in the literature.

 1) To get the topological invariance, it was important and indeed
 necessary for us to impose the
 condition that $mB$ be objective  in Theorem \ref{chern}. For an
 arbitrary flat B-field
$B$, without invoking $L$ or assuming $mB$ to be objective, one
can still prove that $\phi(E;\mc{A},B,D)=\det(I+\frac{1}{2\pi
i}\wt{F}_D)$ is closed. (Indeed one has a Bianchi identity
directly from the following computation
$$d\wt{F}_D=d(F_{D_i}-B_iI)=dF_{D_i}=[F_{D_i},\theta]=[\wt{F}_D,\theta],$$
where $\theta$ is the connection matrix of $D_i$ under a local
frame and the last equation follows from $[B_iI,\theta]=0$ since
$B_i$ is a 2-form.)
 One could then try to define the Chern class of $E$ as
$\wt{c}_k(E)=[\phi_k(E;\mc{A},B,D)]$. In terms of
 bundle gerbe modules, this seems to correspond to the Chern character
using an arbitrary curving $f$ in Section 6.3 of Bouwknegt {\em et
al} \cite{bs}. However such an expanded definition in either place
 has the issue that $\wt{c}_k(E)$ is not well-defined
but depends on the choices made for $\mc{A},B,D$. Without invoking $L$,
it is not possible to characterize such dependence.
The indeterminacy is essentially due to the fact that the compatibility between
$B$ and $\mc{A}$ only mildly constrains $B$ by $\mc{A}$.

2) In \cite{ca}, C$\breve{\mbox{a}}$ld$\breve{\mbox{a}}$raru
proposed to define the Chern classes $c^W(E)$ of $E$ as the
regular Chern classes $c(E\otimes W)$, where $W$ is a fixed
twisted vector bundle over $\mc{L}^{-1}$. The main
 problem in this approach is that no natural choices  were given for $W$ in
applications. In comparison, our choice $L$ as a trivialization of
$\mc{L}^m$ is pertinent and natural, and does indeed reflect the fact that
$\mc{L}$ is an $m$-torsion as the underlying gerbe of $E$.
  Compare also with \cite{hs}, where $B$ is taken to be a global
closed 2-form on $X$. However this works only because the gerbe $\mc{L}$
is trivial for their applications in K3-complex surfaces $X$.

3) When the gerbe connection $\mc{A}$ is flat, namely all
$F_{A_{ij}}=0$, one can take $B_i=0$ for all $i$. (This is a
constraint on $L$.) Thus $\{F_{D_i}\}$ is already a global section
of $\mr{End}(E)$ and our Chern classes $c^L_k(E)$ are given by
$\det(I+\frac{1}{2\pi i}F_{D_i})$. Passing over to bundle gerbe
modules, this special case should correspond to the Chern classes
in \cite{mms1,ms}.
\end{rem}

It is an interesting problem to apply the objective Chern-Weil theory
to various index theorems and twisted K-theory, which we hope to return
in a future work.

%
%

\section{Objective sections and subbundles of twisted bundles}\label{oss}

We begin by the following definition; the second part seems to be brand new.

\begin{df}\label{sec} Suppose $E\prec\mc{L}$ and $L\prec\mc{L}^m$.
As in formula (\ref{m3}), $\mc{L}$ comes with nowhere vanishing
local sections $\xi=\{\xi_{ij}\}$.

1) A {\em twisted section} of $E$ is a collection of local
sections $s=\{s_i\}$ of $E_i$ satisfying $s_j=s_i\otimes\xi_{ij}$
under the identification $E_j=E_i\otimes L_{ij}$.

2) Note $L^{-1}\prec\mc{L}^{-m}$ and the latter carries local
sections $\xi^{-m}=\{\xi_{ij}^{-m}\}$. Suppose $\xi^{-m}$ is
objective, namely there are local sections $t=\{t_i\}$ of $L^{-1}$
such that $t_j=t_i\otimes\xi_{ij}^{-m}$. We will rephrase $s$ as
a  {\em $t$-objective section} of
$E$ in this context, in order to emphasize the relevance of
$\{t_i\}$.
\end{df}

The main point of introducing objective sections is that
$s_i\otimes \sqrt[m]{t_i}$ will be a global section, namely
$s_i\otimes \sqrt[m]{t_i}=s_j\otimes \sqrt[m]{t_j}$ on the
overlap, whenever $\sqrt[m]{t_i}$ is  properly defined. Of course
here $t$ is just a twisted section of $L^{-1}$ over $\mc{L}^{-m}$.
We use $L^{-1}$ instead of $L$, since $t_i$ may vanish so
$t_i^{-1}$ may not exist.

Since generally $\xi^{-m}$ is not objective, the existence of
$t$  does constrain $s$ via $\xi$. In fact one checks readily the
existence of a nowhere vanishing $t$ is equivalent to that the
$\check{\mbox{C}}$ech 2-cycle $\un{z}=\{z_{ijk}\}$ of $\mc{L}$
 has order $m$:  $\un{z}^m=1$. By Proposition
\ref{tor}, such a 2-cycle exists only on special covers.

%
%

From sections we can define twisted bundle homomorphisms in a
rather formal way.

\begin{df}\label{tho}
Suppose $F\prec\mc{K}, E\prec\mc{L}$ are twisted vector bundles of
ranks $n, m$ respectively. Let $K\prec\mc{K}^n, L\prec\mc{L}^m$ be
chosen trivializations.

1) A {\em twisted homomorphism} $f: F\to E$ is just any twisted
section of the bundle $F^*\otimes E\prec\mc{K}^*\otimes\mc{L}$.

2) From a twisted homomorphism $g:L^n\to K^m$, one has then a
{\em $g$-objective homomorphism} $f$.
Here $L^n\prec\mc{L}^{mn}, K^m\prec\mc{K}^{nm}$ are the associated
twisted bundles. (Note that $g$ maps in the opposite direction as
$f$.)
\end{df}

More clearly, the gerbes $ \mc{K}, \mc{L}$ come with local
trivializations $\eta_{ij}, \xi_{ij}$ as in Equation (\ref{m3}).
They give rise to bundle isomorphisms $\chi_{ij}: K_{ij}\to
L_{ij}$ by sending $\eta_{ij}$ to $\xi_{ij}$. Then $f=\{f_i:
F_i\to E_i\}$ should satisfy
$$f_j=f_i\otimes \chi_{ij}\;\;\mr{ on }\;\; F_j=F_i\otimes K_{ij}.$$
Moreover $g=\{g_i: L^n_i\to K^m_i\}$ and $g_j=g_i\otimes
\chi_{ij}^{-mn}$ on $\;\; L_j^n=L_i^n\otimes L_{ij}^{mn}$, where
$\chi_{ij}^{-1}$ is the inverse of $\chi_{ij}$.
  For a $g$-objective homomorphism $f$, $f\otimes\sqrt[mn]{g}$ will be a
global homomorphism if suitably defined.

In particular we can define a twisted or $g$-objective subbundle
$F\prec\mc{K}$ of $E\prec\mc{L}$. However in this approach, it is
unclear how $\mc{K}$ should be tied to $\mc{L}$. Obviously one can
not always choose the same gerbe $\mc{K}=\mc{L}$ for all
subbundles. For example for a twisted line bundle $F$, the gerbe
$\mc{K}$ must be trivial. Thus if $\mc{L}$ is non-trivial, then
$E$ does not allow any twisted line subbundles if one insists on
using the same gerbe.

To handle the issue more decisively, we shall adopt a different
approach for twisted subbundles and make a comparison with the
above general approach at the end. We start by the following
arithmetic result. All integers are assumed to be positive.

\begin{lem}\label{hv}
Given two integers $n\leq m$, let $T=T_{n,m}$ denote the set of
integers $d$ such that each  $nd$ is divisible by $m$. Then $T$
contains two special elements $m_*=\frac{m}{\mr{gcd}(m,n)}$ and
$m^* =C(m,n)$. Moreover $m_*$ is the common divisor hence the
smallest in $T$.
\end{lem}

The last statement is proved by factoring $m,n$ into products of
primes. (Of course $m\in T$ and $m_*\leq m\leq m^*$, hence the
notations.) The key issue in defining twisted subbundles is to
spell out what should be the underlying twisting gerbes.

\begin{df}\label{sb}
Suppose $E=\{E_i\}\prec\mc{L}$ is a rank $m$ twisted vector
bundle. An {\em $\mc{L}^d$-twisted subbundle} $F$ of rank $n$ is a
collection of subbundles $\{F_i\}$ of $\{E_i\}$ such that $F$ by
itself is a twisted bundle over $\mc{L}^d$, where $d\in T_{n,m}$.
 \end{df}

The definition makes sense since each $\mc{L}^d$ is indeed an
$n$-torsion gerbe by the definition of $T$. If $n=1$, then
$m_*=m^*=m$ and any twisted line subbundle in this case is over
the gerbe $\mc{L}^m$, hence a trivialization of $\mc{L}^m$. For
other $d\in T$, an $\mc{L}^d$-twisted line subbundle of $E$ may be
subject to a bigger gerbe $\mc{L}^{ms}$ for some $s$. The
proposition below exhibits certain compatibilities with wedge
product and subbundle operations.

\begin{pro}\label{sswb} Suppose $E\prec\mc{L}$ is a rank $m$ twisted bundle
and $F, W$ are twisted bundles of ranks $n,r$.

 {\em 1)}  If $F\subset E$ is an $\mc{L}^{m^*}$-twisted subbundle, then
the wedge product $\wedge^k F$ is also  an $\mc{S}^{m^*}$-twisted
subbundle of $\wedge^kE$, where $\mc{S}=\mc{L}^k$ is the twisting
gerbe of $\wedge^nE$.

{\em 2)} If $F\subset E, W\subset F$ are respectively $\mc{L}^d$
and $\mc{G}^{d'}$-twisted subbundles, where $d\in T_{n,m}, d'\in
T_{r,n}$ and $\mc{G}=\mc{L}^d$ is the twisting gerbe of $F$, then
$W\subset F$ is a $\mc{L}^{dd'}$-twisted subbundle.
\end{pro}

\nt{\em Proof.} 1) Since $\wedge^k F$ is over the gerbe
$(\mc{L}^{m^*})^k=\mc{S}^{m^*}$, to show $\wedge^k F\subset
\wedge^kE$ is an $\mc{S}^{m^*}$-twisted subbundle, one only needs
to check $m^*\in T_{n', m'}$ where $n'=\mr{rank}(\wedge^k
F),m'=\mr{rank}(\wedge^kE)$. This in turn follows from the
identity:
$$C(n,k)C(m,n)=C(m-k,n-k)C(m,k).$$
(However $m^*\not=(m')^*$, the latter being $C(m',n')$ by
definition.)

2) This amounts to the natural map: $$\begin{array}{ccc}
T_{n,m}\times T_{r,n}&\longrightarrow&T_{r,m}\\
(d,d')&\mapsto&dd'
\end{array}$$
which can be confirmed directly. q.e.d.

Note that statement 1) is false if $\mc{L}^{m^*}$ is replaced by
$\mc{L}^{m_*}$. Namely for an $\mc{L}^{m_*}$-twisted subbundle
$F\subset E$, the wedge product $\wedge^kF\prec{\mc{S}}^{m_*}$
 may not be a twisted subbundle of $\wedge^kE$, as $m_*\notin
T_{n',m'}$ in general. (For instance take $m=6,n=4, k=3$.)

\begin{rem}\label{detm}
 For applications of subbundles, one often fixes a choice of
$d\in T$. Proposition \ref{sswb} indicates that there is an
advantage in selecting $d=m^*$. For example, given  an
$\mc{L}^{m^*}$-twisted subbundle $F\subset E$, its determinant
$\det F$ is then a twisted line bundle over the gerbe
$\mc{S}^{m^*}$. This is consistent with $\det F\subset \wedge^nE$
being an $\mc{S}^{(m')^*}$-twisted line subbundle, where
$(m')^*=m'=\mr{rank}(\wedge^nE)=m^*$.
\end{rem}

One can easily interpret $\mc{L}^d$-twisted subbundles using
twisted homomorphisms of \ref{tho}.

\begin{pro}\label{cam} Suppose $F\subset E$ is an $\mc{L}^d$-twisted
subbundle.

{\em 1)}  Then $F$ is the same as an injective twisted
homomorphism $f: F\to E$, which satisfies $f_j=f_i\otimes\chi_{ij}
$ on $F_j=F_i\otimes L^d_{ij}$, where $\chi_{ij}: L^d_{ij}\to
L_{ij}$ maps the basis $\xi_{ij}^d$ to the basis $\xi_{ij}$.

{\em 2)} Assume the defining cycle $\un{z}$ of $\mc{L}$ is of
order $m$ so that $L_i$ carries a basis $e_i$ with
$e_j=e_i\otimes\xi^d_{ij}$ by Proposition \ref{tor} . Then $F$ is
objective in the sense that $f$ is a $g$-objective homomorphism,
where $g=\{g_i: L^n_i\to L_i^{nd}\}$ and $g_i$ maps the basis
$e^d_i$ to $e_i^{nd}$.

\end{pro}

\nt{\em Proof}. 1) It follows from $F\prec\mc{L}^d$.

2) Since $F\prec\mc{L}^d$, we have $\mc{K}=\mc{L}^d$ and $K=L^k$
in the notations of \ref{tho}, where $nd=mk$. Hence  $g$ here maps
$L^n$ to $K^m=L^{km}=L^{nd}$.  q.e.d.

Note that $f, g$ depend on the choices of the trivializations
(bases) $\{\xi_{ij}, e_i\}$. When $E=\{E_i\}$ is a collection of
trivial bundles, there is an associated collection of bases
$\{e_i\}$ on $L=\{\det E_i\}\prec\mc{L}^m$.

Twisted subbundles are a subtle issue, but have been neglected in
the literature so far. We should find both definitions \ref{tho}
and \ref{sb} useful in different contexts and can feel free to
switch between them accordingly.

\section{Twisted Hermitian-Einstein metrics and
stability}\label{thea}

A Hermitian gerbe metric $\mc{H}=\{h_{ij}\}$ on $\mc{L}=\{L_{ij}\}$ is any
family of fiber metrics $h_{ij}$ on $L_{ij}$ such that all
trivializations $s_{ijk}$ of $L_{ikj}$ have norm 1 under the induced
product metrics $h_{ijk}$.

Take now $(X,\Phi)$ to be a K\"{a}hler manifold of complex
dimension $n$. Fix a holomorphic gerbe $\mc{L}$, namely
all $L_{ij}$ and $s_{ijk}$   are  holomorphic. In this case the
gerbe class $[z_{ijk}]$ lives in $H^2(X,\mc{O}^*)$, the sheaf cohomology of
non-vanishing holomorphic functions. Then
there is a unique connection
$\mc{A}=\{A_{ij}\}$ on $\mc{L}$ such that each $A_{ij}$ is compatible with
 ${H_{ij}}$ and the holomorphic structure on $L_{ij}$ in the usual sense.
Suggested by the curvature of $\mc{A}$, we say further
that a B-field $B=\{B_i\}$ is {\em compatible with $\mc{H}$ and the holomorphic
structure on $\mc{L}$} if each $B_i$ is purely
imaginary and $B_i^{0,2}=0$ respectively. Thus $B_i$ must be
purely imaginary $(1,1)$-forms. However there are infinitely many
compatible B-fields, unlike the case of a unique compatible connection
$\mc{A}$.

Choose a compatible B-field
$B\in\Om^{1,1}\cap\i\Om^2$. Let $E\prec\mc{L}$ be a twisted
vector bundle of rank $m$ with a Hermitian
metric (which restricts to $\mc{H}$ on $\mc{L}$).

\begin{df}\label{hy}
Suppose $D\prec\mc{A}$  is a $(1,1)$-Hermitian connection on $E$,
namely all its curvature $F_{D_i}\in \Om^{1,1}(\mr{End}(E_i))$.
Then $D$ is called {\em $B$-twisted Hermitian-Yang-Mills}, if its curvature
satisfies
\begin{equation}\label{hyk}
\i\Lambda \wt{F}_D=\i\Lambda(F_D-BI)=cI \end{equation} for some
real constant $c$ on $X$. Here  $\Lambda \wt{F}_D=\wt{F}_D\cdot\Phi$
is the usual projection $\Om^{1,1}\to\Om^0$ along the direction of
$\Phi\in\Om^{1,1}$.
\end{df}

Let us do some preliminary analysis of the nature of the equation (\ref{hyk}).
Since $BI$ is diagonal with equal entries,  (\ref{hyk}) splits into a pair of
equations:
\begin{equation}\label{heh1}
\Lambda F_D=\frac{\mr{tr}(\Lambda F_D)}{m}I,
\;\;\i\:\mr{tr}(\Lambda F_D)-\i\ m\Lambda B=mc.
\end{equation}
By taking a refinement cover we may further assume
$D$ to be an $\mr{SU}(m)$-connection. Then (\ref{heh1}) in turn reduces to
 equations
\begin{equation}\label{heh2}
\Lambda F_D=0,\;\;\Lambda B=\i c.
\end{equation}
 In fact only the first equation is essential, since
  $B$ is already fixed at the beginning and now it also
 needs to satisfy the second equation (which just means
 $B_i=\i c \Phi|_{U_i})$.  More plainly, the first equation in
(\ref{heh2}) is  elliptic on each open set $U_i$, despite that
it is  not well-defined on $X$ as the restrictions $\Lambda F_{D_i}=0,
\Lambda F_{D_j}=0$ do not match on $U_i\cap U_j$. While the second
algebraic equation in $B$  plays an auxiliary role only.  It
is introduced to compensate the first equation so together they
compose the globally well-defined equation in (\ref{hyk}). In
other words, (\ref{hyk}) amounts to solving the same kind of local
elliptic differential equations $\Lambda F_{D_i}=0$ as in the original
gauge theory, once we choose one compatible B-field as a solution
of $\Lambda B=\i c$.

 In order to make use of the objective Chern class in Definition
\ref{tc}, we need to introduce  the objective version of the concept: Suppose
$L\prec{\mc{L}}^m$ is a holomorphic trivialization  carrying a compatible
Hermitian connection $A$. (The holomorphic gerbe class $[\mc{L}]\in H^2(X,\mc{O}^*)$
is an $m$-torsion, for example using the holomorphic twisted line bundle $\det E$.)
Set $B=\frac{1}{m}F_A$. Since $B$ consists of $(1,1)$-purely imaginary forms,
it is a B-field compatible with the Hermitian metric and holomorphic structure
on $\mc{L}$. In this special case we say $D$  is
{\em $A$-objective Hermitian-Yang-Mills}.
Thus we have $\Lambda \wt{F}_D=\Lambda(F_D-\frac{1}{m}F_A)=cI$ and $c$  is given
specifically as
\begin{equation}\label{cst}
c=\frac{2n\pi}{n!\mr{vol(X)}}\frac{\deg_L(E)}{m}
\end{equation}
in terms of the $L$-twisted degree $\deg_L(E):=\int
c^L_1(E)\Phi^{n-1}$.

Equivalently, one can switch the point of view by using metrics  in place of
connections:

\begin{df}\label{he} Let $L\prec\mc{L}^m$ be a holomorphic trivialization
and $\mc{H}^m$ be objective by a Hermitian metric $H$ on $L$.
Suppose $E\prec\mc{L}$ is a twisted holomorphic bundle with a Hermitian metric
$h\prec\mc{H}$.
Then $h$ is called  $H$-{\em objective Einstein}
if the twisted curvature  of the unique associated connection of $h$ satisfies
\begin{equation}\label{heh}
\i\Lambda \wt{F}_h=\i\Lambda(F_h-\frac{1}{m}F_H I)=cI.
\end{equation}
Here $F_H$ is the curvature of the $H$-compatible connection on $L$; the Einstein
 constant factor $c$ of $E$ (with respect to $L$) is given by (\ref{cst}).
\end{df}

\begin{rem}
Since we don't assume the associated connection of $h$
 to be $\mr{SU}(m)$, we will not attempt to solve  the two equations in (\ref{heh2})
separately with $B=\frac{1}{m}F_H$.
 Instead, we fix $H$ as above and consider the consequences when the single
equation (\ref{heh}) does have a solution for $h$. As a matter of fact
it will prove quite useful to  combine $h, H$ together and
view them as a single mteric locally in small neighborhood of any point.
\end{rem}

Naturally $E$ is called objective Hermitian-Einstein if it admits an
objective Einstein
metric with respect to some $H$ and $L$. A few basic facts are recorded
here for later use.

\begin{pro}\label{hepro}
{\em 1)} Every twisted line bundle $K$ is objective Hermitian-Einstein.

{\em 2)} If $E, E'$ are objective Hermitian-Einstein with constant factors $c,c'$,
then so are $E^*, \wedge^rE, E\otimes E'$, with factors $-c, rc, c+c'$ respectively.

{\em 3)} Assume $E, E'$ are twisted over the same gerbe. Then  $E\oplus E'$ is
objective Hermitian-Einstein iff $E, E'$ are so  with equal constant factors.
\end{pro}

\nt{\em Proof}. It is straightforward to check the statements. For future reference,
let us just indicate the trivializations used to compute the various constant factors.
For convenience set $\varphi=\frac{2n\pi}{n!\mr{vol(X)}}$, which is a value
depending on the K\"{a}hler metric $\Phi$.

1) Suppose $K\prec\mc{L}$ with a trivialization $L\prec\mc{L}$. Then the
Einstein constant is
$\varphi\:{\deg_L(K)}=\varphi\: c_1(K\ominus L)$.

2) Suppose $E\prec\mc{L}, E'\prec\mc{L}'$ are twisted bundles of ranks $m,m'$.
Choose trivializations $L\prec\mc{L}^m, L'\prec\mc{L}'^{m'}$.
Then we have $E^*\prec\mc{L}^{-1}, \wedge^r E\prec\mc{L}^r,
E\otimes E'\prec\mc{L}\otimes\mc{L}'$,
and their corresponding trivializations are $L^{-1}, L^w, L^{m'}\otimes L'^m$, where
$w$ denotes $C(m-1,r-1)$.
The Einstein constants of $E^*, \wedge^rE, E\otimes E'$ are respectively given as
$$-\varphi\:\deg_L(E)/m, r\varphi\deg_L(E)/m,
\varphi\:[\deg_L(E)/m+\deg_{L'}(E')/m'].$$

3) Here $E\prec\mc{L}, E'\prec\mc{L}$, and also $E\oplus E'\prec\mc{L}$.
However we may have
totally different trivializations $L\prec\mc{L}^m, L'\prec\mc{L}^{m'}$; after
all $E, E'$
may have different ranks $m, m'$ in the first place. Then $E\oplus E'$ inherits the
trivialization $L\otimes L'\prec\mc{L}^{m+m'}$.
The condition on their Einstein constants means that $\deg_L(E)/m=\deg_{L'}(E')/m'$.
Then the Einstein constant of $E\oplus E'$ equals $\varphi$ times this common value.
 q.e.d.
\begin{rem}
Suppose $(E,h)\prec(\mc{L},\mc{H}), (L,H)\prec(\mc{L}^m,\mc{H}^m)$
as above and $h$ is
$H$-objective Einstein. For a second trivialization $(L',H')
\prec(\mc{L}^m,\mc{H}^m)$,
it is simple to show that $h$ is $H'$-objective Einstein iff the difference metric
$\hat{H}=H'\ominus H$ on $\ell=L'\ominus L$ is Einstein in the usual sense.
Since the holomorphic
line bundle $\ell$ always admits an Einstein metric, $E$ being objective
Hermitian-Einstein is
essentially independent of the choice of a trivialization $L\prec\mc{L}^m$
(with a suitable
choice of metrics on $L$).
\end{rem}

Next we consider the algebro-geometric counterpart.

\begin{df}\label{sty}
Suppose $\mc{L}, L, E$ are as in \ref{he}.  Then $E$ is called
{\em $L$-objective stable} (or just $L$-stable) if for every holomorphic
$\mc{L}^{mw}$-twisted line
subbundle $K\subset \wedge^r E $ $, 1\leq r< m$, we have
\begin{equation}\label{stbl}
C(m,r)\deg_{L^w}(K)<\deg_{L^w}(\wedge^rE)
\end{equation}
where $w=C(m-1,r-1)$ as in Corollary \ref{c1}.
\end{df}

The word ``objective'' is employed in the definition because of the involvement of $L$.
To be sure we have the degrees
$$
\deg_{L^w}(K)=\int
c^{L^w}_1(K)\wedge \Phi^{n-1},\; \deg_{L^w}(\wedge^rE)=\int
c^{L^w}_1(\wedge^rE)\wedge \Phi^{n-1}.
$$
Note that the same trivialization $L^w$ is used for both Chern
classes here. By Corollary \ref{c1}, it is possible to express
(\ref{stbl}) as
\begin{equation}\label{tbc}
\frac{\deg_{L^w}(K)}{r}<\frac{\deg_{L}(E)}{m}.
\end{equation}

One can also define the notion of objective semi-stability by
replacing $<$ with $\leq$ in the formulas (\ref{stbl}) or
(\ref{tbc}).

\begin{pro}\label{bac}
Let $E\prec\mc{L}, Q\prec\mc{S}$ be   twisted holomorphic bundles of
ranks $m, 1$ respectively .

{\em 1)} $Q$ is objective stable.

{\em 2)} $E$ is objective stable iff $E^*$ is objective stable.

{\em 3)} $E$ is $L$-objective stable iff $E\otimes Q$ is $L\otimes Q^m$-objective stable.

{\em 4)} Suppose $E'\prec\mc{L}$ is another twisted bundle with the slope
equal to that $E$.
If $E, E'$ are both $L$-semistable, then  $E\oplus E'$ is $L$-semistable.
(But it can not be $L$-stable unless their degrees are zero.)

\end{pro}

\nt{\em Proof}. All statements are generalizations of the standard
stability and can be proved in a similar fashion.

1) This is evident since $\wedge^1 Q=Q$ does not have any nontrivial
proper sub-twisted bundle.

2) As in the usual case, one can use quotient twisted bundles (sheaves)
and a short exact sequence to prove the statement.

3) Note $\wedge^r(E\otimes Q)=\wedge^rE\otimes Q^r$ and any
sub-twisted bundle $K\subset \wedge^rE\otimes Q^r$ can be written
uniquely as $K=K'\otimes Q^r$ for some sub-bundle $K'\subset \wedge^rE$,
since $Q^r$ is a twisted line bundle.

4) Of course $E\oplus E'\prec\mc{L}$.
Here any $K\subset \wedge^r(E\oplus E')=\wedge^rE\oplus\wedge^rE'$ has a
unique decomposition $K''\oplus K'$. The rest of the proof is clear by working
componentwise.
 q.e.d.

The $L$-stability of a bundle does depend on the choice of
the trivialization $L$ up to a point. More precisely we have the following:

\begin{pro}\label{prt}
Suppose $E\prec\mc{L}$ and $L, L'\prec\mc{L}^m$ with difference
bundle $\ell=L\ominus L'$.

{\em 1)} If $E$ is $L$-stable and $\deg(\ell)=\int c_1(\ell)\wedge
\Phi^{n-1}\leq 0$, then $E$ is also $L'$-stable.

{\em 2)} Let $\mc{S}$ be a gerbe given by
$\sqrt[m]{\ell^{-1}}$ (see Example \ref{e1}) and  $Q\prec\mc{S}$ the associated
 trivialization. Then $E$ is
$L$-stable iff $E\otimes Q$ is $L'$-stable.

{\em 3)} When the gerbe (class) $\mc{L}$ is trivial, the objective
stability corresponds to the standard stability.
\end{pro}

\ni{\it Proof.} 1) Let $K\subset \wedge^rE$ be any $\mc{L}^{mw}$-twisted line
subbundle. By Corollary \ref{c1}, $\deg_{L'}(K)=\deg_L(K)+w\deg(\ell)$
and $\deg_{L'}(\wedge^rE)=\deg_L(\wedge^rE)+w\deg(\ell)$. The results
follows easily.

2) One just needs to note that $Q^m$ is identified naturally with $\ell^{-1}$.
Then by Proposition \ref{bac}, $E$ is $L$-stable iff $E\otimes Q$ is
$L\otimes Q^m=L'$-stable.

3) By assumption, there is a trivialization $J\prec\mc{L}^{-1}$.
 Thus  $W:=E_i\otimes J_i$ is a global holomorphic bundle.
  Then
  one shows, via  Proposition \ref{bac}, that $E$ is objective
  stable iff $W$ is stable in the usual sense.
  (The usual stability condition actually involves all proper
  coherent subsheaves $V\subset W$
  of arbitrary ranks $r$, from which $\det V\subset\wedge^rE$. Then
  one needs to apply  the familiar fact
 that any reflexive rank-1 sheaf is locally free, i.e. a line bundle.
 Compare also with the Appendix of \cite{dr}.)
 q.e.d.

\section{Stability of objective Hermitian-Einstein bundles}\label{ghk}

We  establish here one direction of the generalized
Hitchin-Kobayashi correspondence that every objective
Hermitian-Einstein bundle is objective stable. More precisely we
have the following result. We follow the proof of L\"{u}bke
\cite{l} closely as described by Chapter V of Kobayashi \cite{k}.

\begin{thm} Consider a rank $m$-holomorphic twisted bundle  $E\prec\mc{L}$
together with a trivialization $L$ of $\mc{L}^m$.
 If there is a  Hermitian metric $H$ on $L$
such that $E$ admits an $H$-objective Einstein
 metric, then $E$ is $L$-twisted semistable. Moreover
  $E=E^1\oplus E^2\oplus\cdots\oplus E^k$, where $E^1,\cdots,E^k$
  are all twisted bundles over $\mc{L}$, objective stable,
 and objective Hermitian-Einstein with equal constant factors.
\end{thm}

 \nt{\em Proof}. Take any $\mc{L}^{mw}$-twisted
holomorphic line subbudle $K\subset \wedge^r E$ for each $r$ such
that $1\leq r<m$. In view of (\ref{cst}), the semistable version
of (\ref{tbc}) says that to show $E$ is $L$-twisted semistable, we
need to prove
\begin{equation}\label{main}
rc-c'\geq 0 \end{equation} always holds. Here $c, c'$ are
respectively the constant factors of the $H$-objective Einstein
metrics on $E$ and $K$.

The inclusion $K\subset \wedge^r E$ yields a nowhere vanishing
twisted section $s=\{s_i\}$ of $\hat{E}=\wedge^r E\otimes
K^*\prec\hat{\mc{L}}$, where
$\hat{\mc{L}}=\mc{L}^r\otimes\mc{L}^{-mw}$. Set
$\hat{L}\prec\hat{\mc{L}}^{\hat{m}}$ to be the trivialization
induced by $L$, where $\hat{m}=C(m,r)$ is the rank of $\hat{E}$.
By taking refinement cover if necessary, we can assume $s$ to be
objective (see Definition \ref{sec}), so there are nowhere
vanishing local holomorphic sections $\xi=\{\xi_{ij}\}$ of
$\hat{\mc{L}}$ and $t=\{t_i\}$ of $\hat{L}$ such that
\begin{equation}\label{wel}
s\prec\xi, \;\; t\prec\xi^{\hat{m}}.
\end{equation}
(Here we adjust the notation $t$ slightly based on  the fact that the
holomorphic gerbe $\hat{\mc{L}}$ is of order $\hat{m}$.)

Let $\hat{H}$ be the metric on $\hat{L}$, that is induced by $H$.
By Proposition \ref{hepro}, $\hat{E}$ carries an $\hat{H}$-objective
Einstein metric $\hat{h}$ with the constant factor $rc-c'$. Define
a family of local functions $f=\{f_i\}$ on our manifold $X$, where
each
$$f_i=\| s_i\|^2_{\hat{h}}/\sqrt[\hat{m}]{\|t_i\|^2_{\hat{H}}}\:.$$
Note that $\hat{h}\prec\hat{\mc{H}},
\hat{H}\prec\hat{\mc{H}}^{\hat{m}}$, where $\hat{\mc{H}}$ is the
gerbe metric on $\hat{\mc{L}}$. This and relations (\ref{wel})
together yield the key fact that $f$ is a globally well-defined
(smooth) function on $X$. Thus $f$ must have a maximum $q$ on $X$.

Take any point $x_0\in f^{-1}(q)$ and a coordinate neighborhood
$V$. Introduce the functional on $V$,
$$J(f)=\sum_{\alpha,\beta} g^{\alpha\bar{\beta}}\frac{\partial^2 f}
{\partial z^\alpha\partial \bar{z}^\beta},$$ where
$g^{\alpha\bar{\beta}}$ are the components of the metric $\Phi$.
Now choose $V$ so small that $V\subset U_i$ for some $i$ and
$\sqrt[\hat{m}]{t_i}$ exists as a holomorphic section of
$\sqrt[\hat{m}]{\hat{L}_i}$ on $V$. Thus we have a nowhere
vanishing holomorphic section on $V$,
$$\check{s}=s_i\otimes (t_i)^{-1/\hat{m}}\in \Gamma(\hat{E_i}\otimes
\hat{L}_i^{-1/\hat{m}}),$$ and $f=\|\check{s}\|^2_{\check{h}}$,
where $\check{h}=\hat{h}\otimes \hat{H}^{-1/\hat{m}}$. Apply the
standard Weitzenb\"{o}ck formula to the bundle
$\check{E}_i=\hat{E_i} \otimes\hat{L}_i^{-1/\hat{m}}$ over $V$:
\begin{equation}\label{rert}
J(f)=\sum_{\alpha,\beta} g^{\alpha\bar{\beta}}\frac{\partial^2
\|\check{s}\|^2_{\check{h}}} {\partial z^\alpha\partial
\bar{z}^\beta}=\|D_{\check{h}}\check{s}\|^2_{\check{h}}-
\check{R}(\check{s},\check{s}).
\end{equation}
Here $D_{\check{h}}$ is the unique compatible connection
associated to $\check{h}$ (recall $\check{s}$ is holomorphic), and
$\check{R}$ is its mean curvature.  More precisely $\check{R}$ is
the skew symmetric form corresponding to $\i\Lambda
F_{\check{h}}\in\mr{End}(\check{E}_i)$, where $F_{\check{h}}$ is
the curvature of $D_{\check{h}}$ as before.

It is not hard to check that $\check{h}$ is $\check{H}$-objective
Einstein with factor $rc-c'$ by Proposition \ref{hepro}. (Any
metric $H$ on $L$ is certainly $H$-objective Einstein with factor
$0$ as $\deg_L(L)=0$, and $\hat{H}^{-1/\hat{m}}$ is
 some power of $H$ hence also objective Einstein with factor $0$.)
However by definition,
$\check{H}=\hat{H}\otimes(\hat{H}^{-1/\hat{m}})^{\hat{m}}=1$ is
the trivial metric on the trivial twisting trivialization
$\check{L_i}=\hat{L_i}\otimes(\hat{L_i}^{-1/\hat{m}})^{\hat{m}}
={\bf C}$. Thus $\check{h}$ is Einstein in the usual sense with
the same factor $rc-c'$. This means $\i\Lambda
F_{\check{h}}=(rc-c')I$ and (\ref{rert}) becomes here
\begin{equation}\label{hop}
J(f)=\|D_{\check{h}}\check{s}\|^2_{\check{h}}-
(rc-c')\|\check{s}\|^2_{\check{h}}.
\end{equation}

We now return to the proof of (\ref{main}). Suppose this were not
the case but $rc-c'<0$. Then the formula (\ref{hop}) above would
yield that $J(f)\geq0$ on $V$, in fact $J(f)>0$, as $\check{s}$
never vanishes. Since $J(f)$ attains its maximum at the interior
point $x_0\in V$, the Hopf maximum principle says that $f$ must be
a constant on $V$. Hence $J(f)=0$ on $V$, which is a
contradiction.

Hence $E$ is $L$-twisted semistable. Suppose $E$ is however not
quite $L$-objective stable, namely $rc-c'=0$. Then we will seek the
desirable decomposition as stated in the theorem. First
(\ref{hop}) still says
$J(f)=\|D_{\check{h}}\check{s}\|^2_{\check{h}}\geq0$ on $V$ and
the maximum principle again implies $f$ is a constant on $V$,
which means $V\subset f^{-1}(q)$. This is done for each $x_0\in
f^{-1}(q)$, so $f^{-1}(q)$ is open in $X$ hence must be $X$ as it
is certainly closed. In other words, $f$ is actually a constant
function on $X$.

Now that $f$ is a constant function, $J(f)=0$ on $X$. To each
$x\in X$, from (\ref{hop}), one has $D_{\check{h}}\check{s}=0$ on
a sufficiently small neighborhood $V_x\subset U_i$ for some $i$.
Namely $\check{s}$ is $D_{\check{h}}$-flat on $V_x$. Since we can
arrange easily $t_i$ to be $\hat{H}_i$-flat (as $L_i$ is a line
bundle), we see that $s_i$, namely the inclusion $K_i\subset
\wedge^r E_i$, is $D_h$-flat on $V_x$. It follows that
$E=E'_x\oplus E''_x$ for some subbundles of ranks $r, m-r$ on
$V_x$, with $E'_x, E''_x$ both preserving the connection $D_h$.
Namely the metric $h$ splits with respect to the decomposition:
$E'_x\perp E''_x$. As $x\in X$ varies, we obtain a refinement
cover $\{V_x\}$ of $\{U_i\}$ as well as two bundles $E'=\{E'_x\},
E''=\{E''_x\}$, both of which are twisted over the gerbe $\mc{L}$
restricted to the refinement cover. To this end we have the
decomposition $E=E''\oplus E''$ on $X$, in which $E',
E''\prec\mc{L}$ both carry the induced Einstein metric by $h$.
Repeat the same process for $E', E''$. After finite many steps, we
will have to stop. The final decomposition of $E$ is what is
required in the statement of the theorem.  q.e.d.

\begin{cor}
If $E$ is also indecomposable, namely $E\not=E'\oplus E''$ for any
twisted bundles $E', E''\prec\mc{L}$ of positive ranks, then $E$
is $L$-objective stable.
\end{cor}

In fact,  more is true. One can weaken the assumption by requiring
additionally that $E',E''$ are objective stable in the direct sum.

\section{Existence of objective Einstein metrics on stable bundles}\label{oth}

 Now we work on the other direction of the correspondence,
  following the original papers \cite{d1,d2,uy}
 and the expositions \cite{s, k,li, dk}. As in  \cite{s}, we will adapt the approach
 that combines \cite{d1} and \cite{uy} together to our situation.

\begin{thm}\label{don}
Suppose $E\prec\mc{L}$ is a rank-$m$ holomorphic twisted bundle and
$L\prec\mc{L}^m$ is a trivialization.
If $E$ is $L$-objective stable, then $E$ admits an  $H$-objective Einstein metric $h$
for some Hermitian metric $H$ on $L$.
\end{thm}

\nt{\em Proof}. Although we state the theorem for an Einstein metric $h$, it is the
most convenient to combine $h, H$ together and solve the
Einstein equation for the pair $\wt{h}=(h,H)$. Namely  we begin by introducing
the set
$$\mc{M}=\{(h,H):  h, H \mbox{ are compatible
 twisted Hermitian metrics on }E, L\}.$$
Each $h$ leads to a unique gerbe metric $\mc{H}$ on $\mc{L}$ and
the compatibility means that $H\prec\mc{H}^m$. In a small
neighborhood of any point, it will be essential to view
$\wt{h}\in\mc{M}$ as an  ordinary (untwisted) Hermitian metric
 $\wt{h}_i=h_i\otimes{H_i}^{-1/m}$ on the bundle
$\wt{E}_i=E_i\otimes (\sqrt[m]{L_i})^*$ (where-ever existing). In this way the local
picture matches that of the standard case in \cite{d1, d2, uy, s},
and similar local computations, including pointwise
 algebraic operations, can be
carried over to our case without essential changes. The proof can then be repeated
almost verbatim.
We outline below  the main steps and leave the detailed checking to the interested
reader.

\nt{\em{Step 1}}. Define a Donaldson type functional $\mc{D}:\mc{M}\to\R$.

Dropping the positive definiteness requirement for $(h,H)$ in $\mc{M}$, we have the
larger space $\mc{B}$ of pairs of compatible Hermitian forms on $E, L$. Here
$\mc{B}$ is a Banach
space containing $\mc{M}$ as a convex subset. This endows $\mc{M}$ with a
Banach manifold structure  and identifies the tangent space $T_{\wt{h}}\mc{M}=\mc{B}$
naturally at any point $\wt{h}$. Fix some $\wt{k}=(k,K)\in\mc{M}$ as a base point.
Then $Q_1(\wt{h})=\{\log(\det(\wt{k}_i^{-1}\wt{h}_i))\}$ is defined, since $\wt{k}_i, \wt{h}_i$
are genuine metrics pointwise and $\wt{k}_i^{-1}\wt{h}_i\in\mr{End}(\wt{E}_i)$.
In fact $Q_1(\wt{h})$ is a globally well defined function on $X$ by the internal
compatibility
of $\wt{k},\wt{h}$ and $\mr{End}(\wt{E}_i)=\mr{End}(E_i)$ canonically.
Let $\wt{h}(t)=(h(t),H(t))$ be a curve in $\mc{M}$ joining
$\wt{h}$ to $\wt{k}$. Likewise we see that $\wt{h}(t)^{-1}\pa_t\wt{h}(t)\in
 \Om^0(\mr{End}(E))$ is a global section, where $\pa_t\wt{h}(t)=\frac{d}{dt}\wt{h}(t)$.
 The associated connection
 of $\wt{h}_i$ on $\wt{E}_i$ has the curvature $F_{\wt{h}_i}$ equal to
 the twisted curvature $\wt{F}_{h_i}=F_{h_i}-\frac{1}{m}F_{H_i}I$ of $h$. As it has
 used several times before,
  the collection $\{\wt{F}_{h_i}\}$ yield a global section in $\Om^2( \mr{End}(E))$.
  Hence
we have a global section  $F_{\wt{h}}=\{F_{\wt{h}_i}\}$  as well and we can introduce
a 2-form $Q_2(\wt{h})=\i\int^1_0\mr{Tr}[(\wt{h}(t)^{-1}\pa_t\wt{h}(t))
\cdot {F}_{\wt{h}(t)}]dt$
on $X$. Now define the functional
$$\mc{D}(\wt{h})=\int_X[Q_2(\wt{h})-\frac{c}{n}Q_1(\wt{h})\Phi]\wedge
\frac{\Phi^{n-1}}{(n-1)!}$$
with $c$  the (potential) Einstein constant given in (\ref{cst}). To verify
$Q_2(\wt{h})$ is
independent of the paths $\wt{h}(t)$, it is sufficient to check
$$\oint_c \mr{Tr}[(\wt{h}(t)^{-1}\pa_t\wt{h}(t))\cdot {F}_{\wt{h}(t)}]dt\in
\mr{Im}\pa+\mr{Im}\ov{\pa}$$
for any closed path $c=\wt{h}(t)\subset\mc{M}$. The latter depends on the same
kind of local computations as in the standard case \cite{d1} and they can be
transplanted over directly.

\nt{\em{Step 2}}. Establish the main properties of the grading flows for $\mc{D}$.

Given two tangent vectors $v,w\in \mc{B}=T_{\wt{h}}\mc{M}$, $\wt{h}^{-1}v$ and
$\wt{h}^{-1}w$ are
global sections of $\mr{End}(E)$. The Riemann metric on $\mc{M}$ is defined by
the inner product $(v,w)=\int_X\mr{Tr}(\wt{h}^{-1}v\cdot \wt{h}^{-1}w)\Phi^n$.
For a smooth family of $\wt{h}(t)\subset\mc{M}$, one has
$\frac{d}{dt}Q_1(\wt{h}(t))=\mr{Tr}(\wt{h}(t)^{-1}\pa_t\wt{h}(t))$ and
$$ \frac{d}{dt} Q_2(\wt{h}(t))-
\i\mr{Tr}(\wt{h}(t)^{-1}\pa_t\wt{h}(t)\cdot F_{\wt{h}(t)})\in
\mr{Im}\pa+ \mr{Im}\ov{\pa}$$ by the same local computations as in
\cite{d1,k}. Hence one has the variation of the functional
$$\begin{array}{rl}
\frac{d}{dt}\mc{D}(\wt{h}(t))\!=\!\!&\!\!\int_X[\i\mr{Tr}(\wt{h}(t)^{-1}
\pa_t\wt{h}(t)\cdot F_{\wt{h}(t)})-
\frac{c}{n}\mr{Tr}(\wt{h}(t)^{-1}\pa_t\wt{h}(t))\Phi]\wedge
\frac{\Phi^{n-1}}{(n-1)!}\\
=\!\!&\!\!\int_X\mr{Tr}[\i\wt{h}(t)^{-1}\pa_t\wt{h}(t)\cdot F_
{\wt{h}(t)}\wedge\Phi^{n-1}-
\frac{c}{n}\wt{h}(t)^{-1}\pa_t\wt{h}(t)\Phi^n]\frac{1}{(n-1)!}\\
=\!\!&\!\!\int_X\mr{Tr}[\wt{h}(t)^{-1}\pa_t\wt{h}(t)(\frac{\i}{n}\Lambda
F_{\wt{h}(t)}\Phi^n-
\frac{c}{n}\Phi^n)]\frac{1}{(n-1)!}\\
=\!\!&\!\!\int_X\mr{Tr}[\wt{h}(t)^{-1}\pa_t\wt{h}(t)\cdot\wt{h}(t)^{-1}
\wt{h}(t)({\i}\Lambda F_{\wt{h}(t)}-cI)]\frac{\Phi^n}{n!}\\
=\!\!&\!\!(\pa_t\wt{h}(t), \wt{h}(t)({\i}\Lambda F_{\wt{h}(t)}-cI))
\end{array}$$
which means the gradient vector field of $\mc{D}$ is
$\mr{grad}\mc{D}(\wt{h})=\wt{h}({\i} \Lambda F_{\wt{h}}-cI)$. In
view of $\Lambda F_{\wt{h}}=\Lambda \wt{F}_{{h}}$ and (\ref{heh}),
it follows also that $\wt{h}=(h,H)$ is a critical point of
$\mc{D}$ iff $h$ is an $H$-objective Einstein metric.

A downward gradient flow $\wt{h}=\wt{h}(t)$ of $\mc{D}$ is then determined by
 the equation $\pa_t \wt{h}=-\wt{h}({\i}
\Lambda F_{\wt{h}}-cI)$. Or in terms of the global section $s(t)=
\wt{h}(0)^{-1}\wt{h}(t)$ on $\mr{End}(E)$,
\begin{equation}\label{ket}
 \frac{d{s}}{dt}=-\Delta'_0 s-\i s(\Lambda F_{\wt{h}(0)}-cI)+\i
 \Lambda(\ov{\pa}_0{s}\cdot{s}^{-1}\pa_0{s})
 \end{equation}
where the $\Delta'_0, \pa_0, \ov{\pa}_0$ are the
operators associated with the initial metric pair $\wt{h}(0)$.
 This is a non-linear parabolic equation of the same type as in \cite{d1}
and can be solved essentially in the same way. Briefly, by linearizing the equation
and the Fredholm theory one obtains
the short time existence, and by a continuity argument, one has the
long time existence. In the end, together with the maximum principle
one shows that starting at any $\wt{h}(0)\in\mc{M}$,
the evolution equation (\ref{ket}) has a unique solution $\wt{h}(t)$ defined
for all time $t$
together with a uniform bound
\begin{equation}\label{go}
\max_X|\Lambda F_{\wt{h}(t)}|\leq C.
\end{equation}
For the remaining proof we will set $\wt{h}(0)=\wt{k}$ so $\det
(\wt{k}^{-1}\wt{h}(t))=1$ at
$t=0$. By the maximum principle again, we have $\det(\wt{k}^{-1}\wt{h}(t))=1$
identically
along the entire flow.

\nt{\em{Step 3}}. Under the previous assumption prove the following key estimate
\begin{equation}\label{iu}
\mr{max}_X|\log(\wt{k}^{-1}\wt{h}(t))|\leq C_1+C_2\mc{D}(\wt{h}(t))
\end{equation}
for all $t\geq0$, where $C_1, C_2$ are positive constants independent of $t$.

The $L$-stability is used in this step. The estimate (\ref{iu}) substitutes for a
weaker version in \cite{d2} for the case of projective varieties,  whose proof
requires a Mehta-Ramanathan's restriction theorem on stable bundles.
Here we follow  \cite{s} closely, based on analytic results from \cite{uy}. It
 goes very briefly as follows.

Write $u(t)=\log(\wt{k}^{-1}\wt{h}(t))\in\Om^0(\mr{End}(\wt{E}))$;
this is a self-adjoint endomorphism so with all real eigenvalues.
(We  point out that although $\wt{E}=E_i\otimes (\sqrt[m]{L_i})^*$
involves the choice of the root $\sqrt[m]{L_i}$, the eigenvalues
and the various norms of $u$ are well-defined. This note should be
kept in mind for the rest of the argument.) Because of (\ref{go}),
there are constants $A_1,A_2$ independent of $t$ such that
$\max_X{|u|}\leq A_1+A_2\|u\|_{L^1}$.

We will show (\ref{iu}) using proof by contradiction. Suppose it were not true.
Then there are sequences of
$t_i\to\infty$ and constants $B_i\to\infty$ such that
\begin{equation}\label{mtr}
\|u_i\|_{L^1}=\|u(t_i)\|_{L^1}\to\infty,\;\;\|u_i\|_{L^1}\geq B_i\mc{D}(\wt{h}(t_i)).
\end{equation}
Re-normalize $v_i=u_i/\|u_i\|_{L^1}$ so that $\|v_i\|_{L^1}=1$.
Moreover by (\ref{mtr}), $\ov{\pa}v_i$ is bounded in $L^2$. Thus
there is a subsequence $v_i$ weakly convergent to $v_\infty\in
L^2_1$. The limit $v_\infty$ is self-adjoint almost everywhere and
it can be shown that its eigenvalues $\lambda_1,\cdots,\lambda_r$
are constants. Let $\{\gamma\}$ be the set of intervals between
the eigenvalues and for each $\gamma$ choose a function
$p_\gamma:\R\to\R$ such that $p_\gamma(\lambda_i)=1$ for
$\lambda_i<\gamma$ and $p_\gamma(\lambda_i)=0$ otherwise.

Define $\pi_\gamma=p_\gamma(v_\infty)\in L^2_1(\mr{End}(\wt{E}))$. This can be viewed
as an $L^2_1$ subbundle $\wt{S}_\ga$ of $\wt{E}$ through projection in the sense
that $\pi_\gamma^2=
\pi_\gamma$. By the main regularity theorem of \cite{uy}, $\wt{S}_\gamma$ is
actually a
smooth subsheaf of $\wt{E}$. Since the proof is purely local, there is a smooth
subsheaf $S_\gamma$ of $\{E_i\}$ such that $\wt{S}_\gamma=S_\gamma\otimes
\mc{O}(\sqrt[m]{L_i})^*$
locally. Taking the determinant we set
$$K_\gamma=\det S_\gamma=(\wedge^{b}S_\gamma)^{**}$$
which is an $\mc{L}^{mw}$-twisted line subbundle of $\wedge^b E$,
where $b=\mr{rank}\:S_\gamma=\mr{Tr}(\pi_\ga)$ and $w=C(m-1,b-1)$.
We will show for one of the $K_\gamma$,
\begin{equation}\label{tbe}
\frac{\deg_{L^w}(K_\gamma)}{b}=
\frac{\deg_{L^w}(K_\gamma)}{\mr{Tr}(\pi_\ga)}\geq\frac{\deg_{L}(E)}{m}.
\end{equation}
That is to say,
$C(m,b)\deg_{L^w}(K_\gamma)\geq\deg_{L^w}(\wedge^b E)$.
This however violates the $L$-stability of $E$,
 yielding  our expected contradiction.

To show (\ref{tbe}) let $\wt{k}_\ga$ be the restriction of the
base metric $\wt{k}$ to $\wt{S}_\ga$. By direct computations, the
curvature satisfies $\Lambda F_{\wt{k}_\ga}=\pi_\ga\La
F_{\wt{k}}\pi_\ga+\La(\ov{\pa}\pi_\ga\pa\pi_\ga)$ for every $\ga$.
It follows that $\deg_{L^w}(K_\gamma)=\i\int_X\mr{Tr}(\pi_\ga\La
F_{\wt{k}})\Phi^n-\int_X|\ov{\pa}\pi_\ga|^2\Phi^n$. Let $a$ denote
the biggest eigenvalue of $v_\infty$ and $a_\ga$ the width of the
interval $\ga$. Then $v_\infty=a\mr{Id}-\sum a_\ga\pi_\ga$ and a
combination of degrees can be computed:
$$\begin{array}{ll}
&a\deg_L(E)-\sum_\ga a_\ga\deg_{L^w}(K_\gamma)\\
&=\i\int_X\mr{Tr}(v_\infty)\La F_{\wt{k}}\Phi^n+\int_X \sum_\ga a_\ga|
\ov{\pa}(\pi_\ga)|^2\Phi^n\\
&=\i\int_X\mr{Tr}(v_\infty)\La F_{\wt{k}}\Phi^n+\int_X (\sum_\ga
a_\ga (d p_\ga)^2(v_\infty)(\ov{\pa}v_\infty), \ov{\pa}v_\infty)\Phi^n\\
&\leq0
\end{array}$$
where the last inequality comes from Lemma 5.4 of \cite{s}. On the other hand,
$$a\cdot m-\sum\mr{Tr}(\pi_\ga)=\mr{Tr}(v_\infty)=0.$$
Together with the inequality above, this shows that (\ref{tbe}) must hold
for at least one $\ga$.

\nt{\em{Step 4}}. Existence of the objective Einstein metric as a limit metric.

Owing to the strong estimate (\ref{iu}), unlike \cite{d1}, here we can avoid
 Uhlenbeck's theorems on removable singularities and Columbo gauges  (which do not
hold anyway in higher dimensions).

As before, set $s(t)=\wt{k}^{-1}\wt{h}(t)\in\Om^0(\mr{End}(E))$. Then (\ref{iu})
translates into $\mr{max}_X|\log(s)|\leq C_1+C_2\mc{D}(\wt{h}(t))$, which yields
 two consequences:

  (i) $\mc{D}$ is bounded below along $\wt{h}(t)$ .

  (ii) $\|s(t)\|_{L^p}$ is bounded above for each $p$.

\nt To see (ii), just note $\mc{D}$ is decreasing along the downward flow $\wt{h}(t)$,
hence $s(t)$ is bounded in $C^0$- as well as $L^p$-norms.

 From a minimizing sequence of $\mc{D}$ by (i) and the mean value theorem,
 we have  a sequence $t_i\to\infty$ such that
\begin{equation}\label{gto}
\frac{d}{dt}\mc{D}(\wt{h}(t_i))=-\|\Lambda F_{\wt{h}(t_i)}-cI\|^2_{L^2}\to 0.
\end{equation}
Then from the expression of $\mc{D}$ and noting $Q_1(\wt{h}(t))=0$ identically here,
one sees $\ov{\pa}(s(t_i))$
is bounded in $L^p$. Together with (ii), this implies
 that $s(t_i)$ is bounded in $L^p_1$. Choose $p>n$ and
by the Kondrakov compactness of $L^p_1\hookrightarrow C^0$ there is a subsequence
$s(t_i)$ convergent to $s_\infty$ in $C^0$. Then by (\ref{go}) together
with the gradient equation (\ref{ket})
and the maximum principle, it is possible to show that $|\Delta s(t_i)|$ is uniformly
bounded on $X$ in the sequence. (Compare the proof of Lemma 19 in \cite{d1}.)

Thus for any $p$, $\Delta s(t_i)$ is bounded in $L^p$. By the ellipticity of the
Laplace
operator, this and (ii) together imply that $s(t_i)$ is bounded in $L^p_2$. By taking
a subsequence if necessary, $s(t_i)$ weakly converges to $s_\infty$ in $L^p_2$.
In other words, $\wt{h}(t_i)\rightharpoonup \wt{h}_\infty=\wt{k}s_\infty $
in the $L^p_2$-norm
(with respect to $\wt{k}$ as before). Consequently $F_{\wt{h}_\infty}\in L^p$
exists as
a distribution and $F_{\wt{h}(t_i)}\rightharpoonup F_{\wt{h}_\infty}$ in $L^p$.
In view of (\ref{gto}), we now have $\Lambda F_{\wt{h}_\infty}-cI=0$ weakly. The
standard
elliptic regularity guarantees that $\wt{h}_\infty$ must be smooth, and we have the
desired $H_\infty$-objective Einstein metric $h_\infty$, which is read off
$\wt{h}_\infty=(h_\infty, H_\infty)$.
 q.e.d.

%
%

\section{An application to gauge theory: the $\mr{SO}(3)$-moduli space}\label{aap}

We have so far worked with twisted vector bundles. For slight ease in dealing with
structure groups and
associated bundles, it is also useful to introduce twisted principal bundles
that can be briefly laid out as follows. Fix a central extension of Lie groups
$1\to Z\to G\to H\to1$, for example take
\begin{equation}\label{cex}
1\to\mr{U}(1)\to\mr{U}(n)\to\mr{PU}(n)\to1
\end{equation}
including $1\to \mr{U}(1)\to \mr{U}(2)\to \mr{SO}(3)\to1$ at
$n=1$. In general the multiplication maps $m:Z\times G\to G$ and
$m:Z\times Z\to Z$ are clearly group homomorphisms.

Under the group extension, a {\em principal gerbe} $\mc{P}=\{P_{ij}\}$
consists of a collection of principal $Z$-bundles satisfying
the suitable conditions as in vector bundle gerbes. Here a ``tensor product''
$P_{ij}\otimes P_{jk}$ by definition is
the associated principal $Z$-bundle $(P_{ij}\tilde{\times}P_{jk})\times_m Z$ of
the fiber product via the multiplication homomorphism $m:Z\times Z\to Z$. A
{\em twisted principal bundle}
$Q\prec\mc{P}$ consists of a collection of $G$-bundles $Q=\{Q_i\}$ such that
$Q_j=Q_i\otimes P_{ij}$, where the ``tensor product'' is the associated principal
$G$-bundle
$(Q_i\tilde{\times}P_{ij})\times_m G$  via the homomorphism
$m:Z\times G\to G$. Its projection $\mr{P}(Q)$ via $G\to H$ is obviously a
(untwisted)
 principal $H$-bundle. Of course one can also define twisted principal bundle
 of structure group
 $Z$ (rather than $G$) over $\mc{P}$,
because the tensor product makes sense here. Furthermore there is not much additional
difficulty in defining gerbe connections, B-fields or twisted connections as
in the case
of twisted vector bundles.

To get our objective version, we assume $\mc{P}$ to be $n$-torsion, meaning
there is a principal
twisted $Z$-bundle $P\prec\mc{P}^n$, or equivalently, $n[\mc{P}]=
1\in\check{H}^2(X,\un{Z})$ in the sheaf cohomology of $Z$-valued smooth
functions on $X$. (This is automatic under the
extension (\ref{cex}).) Fix $P$ and a twisted connection $A$ on $P$. This
brings a unique gerbe connection $\mc{A}$ on $\mc{P}$ so that $A\prec\mc{A}^n$. Then
we can look for twisted connections $D$ on $Q$ that obey $D\prec\mc{A}$. Using
the curvatures
of $D, A$ together, we can simply copy Section \ref{occt} to define
{\em $P$-objective Chern classes} $c^P_k(Q)$, which are actually independent
of the choices of $A,D$.

We now return to our main application, which will be  based on
twisted principal bundles under the extension (\ref{cex}). Take
$X$ to be a smooth Riemannian 4-manifold and $S\to X$ a principal
$\mr{SO}(3)$-bundle. Historically the moduli space $M_S$ of
anti-self-dual connections on $S$ has played an important role in
the applications of gauge theory. One of the main issues is to do
with the orientability and orientations of the moduli space and
has been settled in Donaldson \cite{d3} for the important case
that the Stiefel-Whitney class $w_2(S)$ has an integer lift
(namely $S$ is spin$^c$). Here our purpose is to handle the
general case of an arbitrary $w_2(S)$.

More precisely let $\mf{B}_S$ be space of connections on $S$ modulo the gauge
group. Each connection
$\nabla$ on $S$ induces a connection $d_\nabla$ on the adjoint bundle
$\mf{g}_{_S}=S\times_{ad}\mf{so}(3)$. In turn we have the Fredholm operators
$$\delta_\na=-d^*_\na+d^+_\na:\Om^1(\mf{g}_{_S})\to(\Om^0\oplus\Om^2_+)
(\mf{g}_{_S})$$
parameterized by $[\na]\in\mf{B}_S$. Denote by $\La_S\to\mf{B}_S $ the
determinant line
bundle of the family. Then the orientability and orientation of $\La_S$
correspond exactly to
those of the moduli space $M_S$, whenever the latter is smooth.

\begin{pro}\label{otk}
For any $\mr{SO}(3)$-bundle $S$, the associated line bundle
$\La_S$ is always orientable, namely trivial.
\end{pro}

\nt{\em Proof}. This is mainly a  modification of the proof in
\cite{d3} for our twisted bundle set up. There are three
ingredients. First, because of the exact sequence (\ref{cex}), by
lifting $S$ locally one has a twisted principal $\mr{U}(2)$-bundle
$Q=\{Q_i\}$ over a principal gerbe $\mc{P}$. (The projection
$\mr{P}(Q)$ equals $S$ and the original case with integral
$w_2(S)$ corresponds to the trivial gerbe $\mc{P}$.) Choose a
trivialization $P\prec{\mc{P}}^2$ and let
$$\mf{B}^P_Q=\{(D,A): D, A \mbox{ are compatible connections on } Q, P\}/\sim$$
modulo the twisted gauge transformation pairs.
Here the compatibility means that $A\prec\mc{A}^2, D\prec\mc{A}$
for some unspecified gerbe connection $\mc{A}$ on $\mc{P}$.
When a connection $A^\star$ on $P$ is also chosen, each connection
$\na$ on $S$ lifts to a unique twisted connection $D\prec\mc{A}$ on
$Q$ in view of (\ref{cex}). That is to say, one has an injective map
$f=f_{A^\star}:\mf{B}_S\to\mf{B}^P_Q, [\na]\mapsto[D]$, with the image
set $\mf{B}^{A^\star}_Q=\{[D]: (D,A_0)\in\mf{B}^P_Q\}$.

For each pair $\wt{D}=(D,A)\in \mf{B}^P_Q$,  consider locally $\wt{D}$ as
the connection
$D_i\otimes(\sqrt{A_i})^{-1}$ on $Q_i{\otimes}(\sqrt{P_i})^{-1}$ on an open
set $U_i$. In
turn it induces a {\em global} connection
$$d_{\wt{D}}:\Om^0(\mf{g}^P_Q)\to \Om^1(\mf{g}^P_Q)$$
where $\mf{g}^P_Q$ is the associated untwisted vector bundle $(Q{\otimes}
(\sqrt{P})^{-1})\times_\xi\mf{su}(2)$ and $\xi:\mr{U}(2)\to\mf{su}(2)$ is the
adjoint representation composed with the projection. (Note the same can not
be said for $d_D$ without incorporating $A$.) As a matter of fact, one can
identifies $\mf{g}^P_Q=\mf{g}_S$
canonically and $d_{\wt{D}}=d_\na$ at $\wt{D}=(D,A^\star)$. In other words,
$\La_S=f^*\La^P_Q$ is the pull-back of the determinant bundle  of the Fredholm family
$$\delta_{\wt{D}}=-d^*_{\wt{D}}+d^+_{\wt{D}}:\Om^1(\mf{g}^P_{Q})
\to(\Om^0\oplus\Om^2_+)(\mf{g}^P_Q).$$
Thus to show the theorem, it is sufficient to prove $\La^P_Q\to \mf{B}^P_Q$ is trivial.
Since the choice of a trivialization $P\prec\mc{P}^2$ is immaterial to the
discussion above,
 we may as well select $P_o=\det Q=Q\times_{\det}\mr{U}(1)$ from now on.
Then it remains to show that $\La^{P_o}_Q\to \mf{B}^{P_o}_Q$ is trivial.

Next using the homomorphism $\la:\mr{U}(2)\to\mr{SU}(3)$, $u\mapsto\mr{diag}
(u,\det u^{-1})$,
one introduces the twisted principal $\mr{SU}(3)$-bundle
$$Q^+=(Q{\otimes}P^*_o)\times_\la \mr{SU}(3).$$
(This resembles the stablization $E \leadsto E\oplus\det E^*$ from
a $\mr{U}(2)$-vector bundle to an $\mr{SU}(3)$-bundle in the
original proof of \cite{d3}.) The canonical trivialization $\det
Q^+$ is a global line bundle. It follows that $\mf{B}^{P_o}_Q$ can
be identified with the set $\mf{B}_{Q^+}$ of twisted connections
on $Q^+$. Now that $Q^+$ has the structure group $\mr{SU}(3)$,  by
making suitable homotopy computations similar to \cite{d3} and
\cite{dk} (especially 5.4),  one sees that elements in
$H_1(\mf{B}_{Q^+})$ all come from
$[X,\mr{SU}(3)]=K^{-1}(X)/H^1(X)=H^3(X)$. More directly the slant
product over the twisted Chern class $c_2^{P_o}(Q^+)$ gives such a
homomorphism  $H^3(X)\to H_1(\mf{B}_{Q^+})$. Thus each loop $\phi$
in $\mf{B}^{P_o}_Q$ has the form $\phi_\ga$ for some loop $\ga$ in
$X$ via the Poinc\'{a}re duality $H^3(X)\to H_1(X)$. And to show
$\La^{P_o}_Q$ is trivial, one must show the restriction
$\La^{P_o}_Q|_{\phi_\ga}$ is so for any loop $\ga$.

The last ingredient is to use the gluing to exhibit the loop
$\phi_\ga$ of connections in $\mf{B}^{P_o}_Q$. This is a slight
generalization from the case of standard $\mr{U}(2)$-connections
in \cite{d4,d3} to our twisted connection pairs.
 So let $V\to S^4$ be the negative spinor bundle
on the 4-sphere and $J_\la$ be an instanton on $V$, flattened of a
small scale $\la$
around $\infty\in S^4$. Write the pair $\wt{J}_\la=(J_\la,\det J_\la)
\in\mf{B}^{\det V}_V$,
where $V$ is viewed as a twisted bundle.

Choose a twisted bundle $Q'\prec\mc{P}$ with $c_2^{\det
Q'}(Q')=c_2^{P_o}(Q)-1$. Fix a twisted connection $D'$ on $Q'$.
Since $\mf{g}_S, \La^2_+$ are trivial $\mr{SO}(3)$-bundles over
the loop $\ga$, we can choose a lifting isomorphism $\rho:
\mf{g}_S|_\ga \to \La^2_+|_\ga$ as our gluing parameter. Then for
each pair $(D',A')\in\mf{B}^{\det Q'}_{Q'}$ and any point
$x\in\ga$, we have the glued pair
$[(D',A')\#_{\rho(x)}\wt{J}_\la]\in\mf{B}^{P_o}_Q$. Strictly
speaking, it is better to work with the associated twisted vector
bundles of $Q',Q$ and $D',A'$ should be flattened around $x$. In
fact there is not much difference in gluing the twisted
connections, since the gluing is done in a small neighborhood of
$x$ and essentially it is to glue component connections of $D',
A'$ with $J_\la, \det J_\la$. (As usual one should lift $\rho$ to
gluings between $Q'$ and $V$. But the gauge classes
$[(D',A')\#_{\rho(x)}\wt{J}_\la]$ are independent of the
liftings.) Fix a pair $(D',A')$ and  a uniform scale $\la$ over
the compact set $\ga$. The result is a continuous family
$[(D',A')\#_{\rho}\wt{J}_\la]$ in $\mf{B}^{P_o}_Q$, representing
our loop $\phi_\ga$. The main estimates in 3(d) of \cite{d3} gives
a continuous fiberwise isomorphism $j$ in the following diagram:
$$\begin{array}{ccc}
\pi^*\La^{\det Q'}_{Q'}& \xrightarrow{j}&\La^{P_o}_Q\\
\downarrow&&\downarrow\\
\ga\times\mf{B}^{\det Q'}_{Q'}&\xrightarrow{\#_\rho\wt{J}_\la}&\mf{B}^{P_o}_Q
\end{array}$$
where $\pi$ is the obvious projection. Consequently the  bundle restriction
$\La^{P_o}_Q|_{\phi_\ga}$ is isomorphic to
$\pi^*\La^{\det Q'}_{Q'}|_{\ga\times\{[D',A']\}}$,
which is  the pull back of the single fiber $\La^{\det Q'}_{Q'}$
over the chosen point
$[D',A']$ hence trivial. The triviality of $\La^{P_o}_Q|_{\phi_\ga}$ is
established as required.
 q.e.d.

Furthermore, one can settle the orientations by using the excision
to reduce to the case of a K\"{a}hler manifold $X$, where one can use the
natural complex orientation.

\begin{cor}\label{lst}
For each $P$ the moduli space $M_S$ inherits a unique orientation
$\mf{o}_P$. Two such orientations $\mf{o}_{P}, \mf{o}_{{P'}}$ are
the same iff $c^2_1(\ell)$ is even, where $c_1(\ell)\in H^2(X,\Z)$
is the Chern class of the difference line bundle $\ell=P'\ominus
P$.
\end{cor}

Regardless the $\mr{SO}(3)$-vector bundle $S$, it seems to be an interesting problem
to investigate directly the instanton moduli space $M^L(E)$  on a twisted $\mr{U}(2)$-bundle
$E\prec\mc{L}$ with a fixed trivialization $L\prec{L}^2$. Perhaps one could further
 extract Doanldson type invariants. In terms of algebraic geometry, this
suggests the study of the moduli space of $L$-objective stable bundles.

\end{document}